\input amstex
\documentstyle{amsppt}
\magnification=1200
\hoffset=-0.5pc
\nologo
\vsize=57.2truepc
\hsize=38.5truepc

\spaceskip=.5em plus.25em minus.20em

\define\Bobb{\Bbb}
\define\fra{\frak}
\define\HG{G}
\define\hg{g}
\define\GHH{H}
\define\KK{H}
\define\KKK{H}
\define\kk{h}

\define\GH{H}

\define\atiybern{1}
\define\farakora{2}
\define\freudone{3}
\define\hartsboo{4}
\define\howeone{5}
\define\poiscoho{6}
\define\souriau{7}
\define\oberwork{8}
\define\kaehler{9}
\define\kaehredu{10}
\define\lradq{11}
\define\jacobstw{12}
\define\kostasix{13}
\define\lanmaone{14}
\define\mccrione{15}
\define\satakboo{16}
\define\sekiguch{17}
\define\severone{18}
\define\zakone{19}
\topmatter
\title Severi varieties and holomorphic nilpotent orbits
\endtitle
\author Johannes Huebschmann
\endauthor
\affil 
Universit\'e des Sciences et Technologies
de Lille
\\
U. F. R. de Math\'ematiques
\\
CNRS-UMR 8524
\\
F-59 655 VILLENEUVE D'ASCQ, France
\\
Johannes.Huebschmann\@math.univ-lille1.fr
\endaffil
\date{January 28, 2004}
\enddate
\abstract{
Each of the
four critical Severi varieties arises from  
a minimal holomorphic nilpotent orbit in a simple regular rank 3
hermitian Lie algebra and each such variety
lies as singular locus in a cubic---the chordal variety---in 
the corresponding complex projective space;
the cubic and projective space are identified
in terms of holomorphic nilpotent orbits.
The projective space acquires an exotic K\"ahler structure
with three strata, the cubic is an example of an exotic projective variety
with two strata, and the corresponding Severi variety is the closed stratum in
the exotic variety as well as in the exotic projective space.
In the standard cases, these varieties arise also via K\"ahler reduction.
An interpretation in terms of constrained mechanical systems is included.}
\endabstract

\address{\smallskip
\noindent
USTL, UFR de Math\'ematiques, CNRS-UMR 8524
\newline\noindent
F-59 655 Villeneuve d'Ascq C\'edex,
France
\newline\noindent
Johannes.Huebschmann\@math.univ-lille1.fr}
\endaddress
\endkeywords

\endtopmatter
\document
\leftheadtext{Johannes Huebschmann}
\rightheadtext{Severi varieties and holomorphic nilpotent orbits}

\beginsection 1. Introduction

Let $m \geq 2$. A {\it Severi variety\/} is a non-singular  variety
$X$ in complex projective $m$-space $\bold P^m \Bobb C$
having the property that for some point $O \not \in X$
the projection from $X$ to $\bold P^{m-1} \Bobb C$ is a closed immersion,
cf. \cite\hartsboo\ (Ex.~3.11 p. 316),\,
\,\cite\lanmaone,\,\cite\zakone.
Let $X$ be a Severi variety and $n$ the dimension of $X$.
The critical cases are when $m = \frac 32 n + 2$. Zak \cite\zakone\ 
proved that only the following {\it four\/} critical cases occur: 
\newline\noindent
(1.1) $X=\bold P^2 \Bobb C \subseteq \bold P^5 \Bobb C$ (Veronese embedding)
\newline\noindent
(1.2) $X=\bold P^2 \Bobb C \times \bold P^2 \Bobb C \subseteq \bold P^8 \Bobb C$
(Segre embedding)
\newline\noindent
(1.3) $X=\roman G_2(\Bobb C^6) 
= \roman U(6)\big /(\roman U(2)\times \roman U(4))
\subseteq \bold P^{14} \Bobb C$
(Pl\"ucker embedding)
\newline\noindent
(1.4) $X=\roman{Ad}(\fra e_{6(-78)})
\big/(\roman{SO}(10,\Bobb R) \cdot \roman{SO}(2,\Bobb R))
\subseteq \bold P^{26} \Bobb C$.
\newline\noindent
These varieties arise from the projective planes 
over the {\it four\/} real normed division algebras 
(reals, complex numbers, quaternions, octonions)
by complexification, cf. e.~g. \cite\atiybern.

The purpose of the present paper is to exhibit another
interesting geometric feature of the critical
Severi varieties:
There are exactly {\it four\/} simple regular rank 3 hermitian Lie algebras 
over the reals; these result from the  
euclidean Jordan algebras of hermitian ($3\times 3$)-matrices 
over the {\it four\/} real normed division algebras
by the superstructure construction; and {\it each critical Severi variety 
arises from the minimal holomorphic nilpotent orbit 
in such a Lie algebra.\/} 
This provides geometric insight into these varieties
in the realm of singular Poisson-K\"ahler geometry, which
we formulate as Theorem 1.5 below. 

Recall that an {\it exotic K\"ahler structure\/} is a complex analytic 
stratified K\"ahler structure (a definition of the latter will be
reproduced in Section 5 below) with at least two strata, cf. \cite\kaehler. 

\proclaim{Theorem 1.5}
For $m=5,8,14,26$, the complex projective space 
$\Bobb P^m \Bobb C$ carries an exotic normal K\"ahler structure
with the following properties:
\newline\noindent
{\rm (1)}
The closures of the strata constitute an ascending sequence
$$
Q_1 \subseteq 
Q_2 \subseteq  Q_3=\Bobb P^m \Bobb C
\tag1.6
$$
of normal K\"ahler spaces where, complex algebraically, 
$Q_1$ is a Severi variety and
$Q_2$ a projective cubic hypersurface,
the chordal variety of $Q_1$.
\newline\noindent
{\rm (2)}
The singular locus of $Q_3$, in the sense of stratified
K\"ahler spaces, is the hypersurface $Q_2$, and that of $Q_2$
(still in the sense of stratified K\"ahler spaces)
is the non-singular variety $Q_1$; furthermore,  $Q_1$ is as well the
complex algebraic singular locus of $Q_2$.
\newline\noindent
{\rm (3)}
The exotic  K\"ahler structure
on $\Bobb P^m \Bobb C$ restricts to an ordinary K\"ahler structure
on  $Q_1$ inducing, perhaps up to rescaling,
the standard hermitian symmetric space structure. 
\endproclaim

The {\it chordal variety\/} of a non-singular projective variety $Q$
in $\Bobb P^m \Bobb C$ is the closed subspace of $\Bobb P^m \Bobb C$
generated by all chords.
When the chordal variety of a non-singular projective variety $Q$
in $\Bobb P^5 \Bobb C$ is a hypersurface (and not the entire ambient space)
the projection from a generic point gives an embedding in
$\Bobb P^4 \Bobb C$. A classical result of Severi \cite\severone\ 
says that the Veronese surface is the only surface
(not contained in a hyperplane) in
$\Bobb P^5 \Bobb C$ with this property.
This is the origin of the terminology \lq\lq Severi variety\rq\rq.
In \cite\atiybern, the chordal varieties are written as $Z_n(C)$ ($n=0,1,2,3$).

For a compact K\"ahler manifold $N$ which is complex analytically a projective 
variety, with reference to the Fubini-Study metric, the Kodaira embedding 
will not in general be symplectic. The above theorem shows that there are 
interesting situations where {\sl complex  projective space carries an 
exotic K\"ahler structure\/} which, via the {\sl Kodaira embedding, 
restricts to the K\"ahler structure on\/} $N$. More generally,
for a projective variety $N$ with singularities,
the correct question to ask is whether,
via the Kodaira embedding, $N$ inherits,
from a suitable exotic K\"ahler structure on complex projective space,
a (stratified) K\"ahler structure.
In  \cite\kaehredu\ we have developped a K\"ahler quantization scheme
for not necessarily smooth stratified K\"ahler spaces,
including examples of K\"ahler quantization on 
projective varieties with stratified K\"ahler structure.
This procedure applies to the circumstances of the above theorem.
Suffice it to note that ignoring the lower strata means working on
a non-compact space, and this then leads to results which are inconsistent
in the sense that the principle that quantization commutes with reduction is
violated. See \cite\kaehredu\ (4.12) for details.

A complex analytic stratified K\"ahler structure 
on a stratified space $Y$ 
involves among other things a {\it real\/} Poisson algebra
of continuous functions on $Y$ which, on each stratum, restricts to an
ordinary smooth symplectic Poisson algebra.
This Poisson structure is {\it independent\/} of the complex analytic 
structure. An interesting feature of the situation isolated in
Theorem 1.5 is that
the real Poisson algebra on $\Bobb P^m \Bobb C$
(beware: it contains more functions than just the ordinary smooth ones)
detects the lower strata $Q_s \setminus Q_{s-1}$ ($1 \leq s \leq 3$ where
$Q_0$ is understood to be empty)
by means of the rank of the Poisson structure,
independently of the complex analytic structures.
This notion of rank can be made precise by means of the Lie-Rinehart
structure related with a general not necessarily smooth Poisson structure
which we introduced in \cite\poiscoho.

More geometric consequences will be explained later.
In particular, in Section 6 below we will 
show how the ascending sequence (1.6) 
for each of the three classical cases 
results from K\"ahler reduction; this includes,
in particular, a G(eometric)I(nvariant)T(heory) quotient-construction.
When the GIT-quotient is related with the corresponding symplectic quotient
via the standard Kempf-Ness procedure, from Theorem 1.5, we then
obtain geometric insight into the singular structure of
the corresponding symplectic quotients.
This, in turn, leads to an interpretation 
of the constituents of the ascending sequence (1.6)
in terms of constrained systems in mechanics, see Remark 6.6 below.
Here the complex analytic structure does not seem to have a direct 
mechanical significance; it helps understanding the kinematical
description. For issues related with quantization,
the complex analytic structure has its intrinsic significance, though,
cf. \cite\kaehredu.

I am indebted to F. Hirzebruch for having introduced me into
Severi varieties.

\medskip\noindent {\bf 2. Lie algebras of hermitian type}
\smallskip\noindent
Following \cite\satakboo\ (p.~54), we define a (semisimple) Lie algebra
of {\it hermitian type\/} to be a pair $(\fra g, z)$ which consists of a 
real semisimple Lie algebra $\fra g$ with a Cartan decomposition 
$\fra g = \fra k \oplus \fra p$ and a central element $z$ of $\fra k$,
referred to as an $H$-{\it element\/}, such that
$J_z = \roman{ad}(z)\big |_{\fra p}$ is a (necessarily $K$-invariant)
complex structure on $\fra p$. Slightly more generally,
a {\it reductive Lie algebra of hermitian type\/} is a reductive Lie algebra 
$\fra g$ together with an element $z \in \fra g$
whose constituent $z'$ (say) in the semisimple part 
$[\fra g,\fra g]$ of $\fra g$ is an $H$-element for $[\fra g,\fra g]$
\cite\satakboo\ (p.~92). Below we will sometimes refer to $\fra g$ alone
(without explicit choice of $H$-element $z$) as a {\it hermitian\/} Lie 
algebra. For a real semisimple Lie algebra $\fra g$, with Cartan 
decomposition $\fra g = \fra k \oplus \fra p$, we write $G$ for an 
appropriate Lie group having $\fra g = \roman {Lie}(G)$ 
(matrix realization or adjoint realization; both will do)
and $K$ for the (compact) connected subgroup of $G$ with \
$\roman{Lie}(K) = \fra k$; the requirement that $(\fra g,z)$ be of 
hermitian type is equivalent to $G\big / K$ being a (non-compact)
hermitian symmetric space with complex structure induced by $z$.

A real semisimple hermitian Lie algebra $\fra g$ decomposes as 
$\fra g = \fra g_0 \oplus \fra g_1 \oplus \ldots  \oplus\fra g_k$
where $\fra g_0$ is the maximal compact semisimple ideal and where
$\fra g_1,\dots,  \fra g_k$ are non-compact and simple.
For a non-compact simple Lie algebra with Cartan decomposition
$\fra g = \fra k \oplus \fra p$, the  $\fra k$-action on $\fra p$
coming from the adjoint representation of $\fra g$
is faithful and irreducible whence the center of $\fra k$
then is at most one-dimensional; indeed $\fra g$ has an $H$-element 
turning it into a Lie algebra of hermitian type if and only if
the center of $\fra k$ has dimension one. In view of E.~Cartan's
infinitesimal classification of irreducible hermitian symmetric spaces,
the Lie algebras $\fra{su}(p,q)$ ($A_n, n \geq 1$, where $n+1 = p+q$),
$\fra{so}(2,2n-1)$ ($B_n, n \geq 2 $),
$\fra{sp}(n,\Bobb R)$ ($C_n, n \geq 2$),
$\fra{so}(2,2n-2)$ ($D_{n,1}, n>2$),
$\fra{so}^*(2n)$ ($D_{n,2}, n >2$)
together with the real forms
$\fra{e}_{6(-14)}$ and $\fra{e}_{7(-25)}$
of $\fra{e}_{6}$ and $\fra{e}_{7}$, respectively,
constitute a complete list of simple hermitian Lie algebras.

We refer to $(\fra g,z)$ as {\it regular\/}
when the relative root system is of type $\roman C_r$;
see \cite\kaehler\ (Proposition 3.3.2) for details.
Thus $\fra {sp}(1,\Bobb R) \cong\fra {su}(1,1)\cong \fra {so}(2,1)$
and 
$\fra{so}^*(4)$
are the only regular rank 1 simple hermitian Lie algebras;
$\fra {sp}(2,\Bobb R)$, $\fra {su}(2,2)$,
$\fra{so}^*(8)$,
$\fra {so}(2,q)$ ($q \geq 3$)
are the only regular rank 2 simple hermitian Lie algebras;
and
$\fra {sp}(3,\Bobb R)$, $\fra {su}(3,3)$,
$\fra{so}^*(12)$,
$\fra{e}_{7(-25)}$
are the only regular rank 3 simple hermitian Lie algebras.

\medskip\noindent{\bf 3. Jordan algebras}
\smallskip\noindent
As usual, denote by $\Bobb O$ the octonions or Cayley numbers.
For $\Bobb K = \Bobb R,\,\Bobb C,\, \Bobb H,\ \Bobb O$
and $n \geq 1$,
consider the {\it euclidean Jordan algebra\/} 
of hermitian $(n \times n)$-matrices
$\Cal H_n(\Bobb K)$ over $\Bobb K$
($\Cal H_n(\Bobb R) = \roman S^2_{\Bobb R}[\Bobb R^n]$),
with Jordan product 
$\,\circ\,$ given by
$x \circ y = \frac 12 (xy + yx)$
($x,y \in \Cal H_n(\Bobb K)$)
where $n \leq 3$ when 
$\Bobb K = \Bobb O$.
See \cite{\farakora,\,\mccrione,\,\satakboo}
for notation and details.
Here
$\Cal H_1(\Bobb R)\cong\Cal H_1(\Bobb C)
\cong\Cal H_1(\Bobb H)\cong\Cal H_1(\Bobb O)\cong \Bobb R$,
$\Cal H_2(\Bobb O)$ is isomorphic to
the euclidean Jordan algebra $J(1,9)$
arising from the Lorentz form of type $(1,9)$ on $\Bobb R^{10}$,
and $\Cal H_3(\Bobb O)$ is the real exceptional rank 3 Jordan
algebra of dimension 27, also called {\it real Albert algebra\/}.

Any regular simple hermitian Lie algebra
arises from a simple real euclidean Jordan algebra
by the superstructure construction;
we list the cases of interest to us as follows:
$$
\alignedat 3
\Bobb K &= \Bobb R \colon
\quad\Cal H_n(\Bobb R) \oplus \fra {gl}(n,\Bobb R) \oplus \Cal H_n(\Bobb R)
&&= \fra{sp}(n,\Bobb R), 
\quad &&r=n
\\
\Bobb K &= \Bobb C \colon
\quad
{\Cal H_n(\Bobb C)} 
\oplus \fra {sl}(n,\Bobb C) \oplus \Bobb R\oplus \Cal H_n(\Bobb C)
&&= \fra{su}(n,n),
&&r=n
\\
\Bobb K &= \Bobb H \colon
\quad
{\Cal H_n(\Bobb H)} 
\oplus \fra {gl}(n,\Bobb H) \oplus \Cal H_n(\Bobb H)
&&= \fra{so}^*(4n),
&&r=n
\\
\Bobb K &= \Bobb O \colon
\quad
{\Cal H_3(\Bobb O)} 
\oplus \Bobb R h^3 \oplus \fra{e}_{6(-26)} \oplus \Cal H_3(\Bobb O)
&&= \fra{e}_{7(-25)},
&&r=3
\endalignedat
\tag3.1
$$
The correspondence between regular simple hermitian Lie algebras
and simple real euclidean Jordan algebras is actually bijective;
adding to the list (3.1) the cases $\fra {so}(2,q)$ arising from the real 
euclidean Jordan algebras $J(1,q-1)$ we obtain a complete list
of regular simple hermitian Lie algebras but this is not relevant below. 

We now restrict attention to the rank 3 case. Complexification of the 
simple rank 3 euclidean Jordan algebras yields the complex simple rank 
3 Jordan algebras $\roman S_{\Bobb C}^2[\Bobb C^3]$,
$\roman M_{3,3}(\Bobb C)$,
$\Lambda^2[\Bobb C^6]$,
$\Cal H_3(\Bobb O_{\Bobb C})=\Cal H_3(\Bobb O)\otimes \Bobb C$, 
the latter being the complex Albert algebra,
and it is well known that these exhaust the
complex simple rank 3 Jordan algebras, see e.~g. \cite\farakora.
As a side remark we note that the rank 3 tube domains are realized in these 
complex Jordan algebras. In particular, the four regular simple rank 3 
hermitian Lie algebras arise in this fashion from the real rank 3 euclidean 
Jordan algebras, and their Cartan decompositions 
$\fra g = \fra k \oplus \fra p$
have the following form where the decomposition 
(C.2) is spelled out for the reductive hermitian Lie algebra
$\fra u(3,3)$ instead of its simple brother $\fra {su}(3,3)$:
\newline\noindent
(C.1) $\fra{sp}(3,\Bobb R)= \fra u(3) \oplus \roman S_{\Bobb C}^2[\Bobb C^3]$
\newline\noindent
(C.2) $\fra u(3,3)=(\fra u(3)\oplus\fra u(3)) \oplus \roman M_{3,3}(\Bobb C)$ 
\newline\noindent
(C.3) $\fra {so}^*(12)= \fra u(6) \oplus \Lambda^2[\Bobb C^6]$
\newline\noindent
(C.4) $\fra e_{7(-25)}
=(\fra e_{6(-78)}\oplus\Bobb R)\oplus\Cal H_3(\Bobb O_{\Bobb C})$,
$\fra e_{6(-78)}$ being the compact form of $\fra e_6$.
\newline\noindent
The resulting (unitary) $K$-representations on the complex vector spaces
$\fra p$ extend to $K^{\Bobb C}$-representations on $\fra p$,  and these
have the following forms:
\newline\noindent
(R.1) The symmetric 
square of the standard $\roman {GL}(3,\Bobb C)$-representation on
$\Bobb C^3$;
\newline\noindent
(R.2) the standard representation 
of $\roman {GL}(3,\Bobb C) \times \roman {GL}(3,\Bobb C)$ on
$\roman M_{3,3}(\Bobb C)$ via multiplication, from the right and left by 
elements of $\roman {GL}(3,\Bobb C)$, of $(3\times 3)$-matrices;
\newline\noindent
(R.3) the exterior square of the standard 
$\roman {GL}(6,\Bobb C)$-representation on $\Bobb C^6$;
\newline\noindent
(R.4) the classical representation of
$\roman E_6(\Bobb C)$ (extended by a central copy of $\Bobb C^*$)
of complex dimension 27, studied already by E. Cartan.

\medskip\noindent
{\bf 4. Jordan rank and Severi varieties}
\smallskip\noindent
Let $\fra g$ still be a simple rank 3 hermitian Lie algebra, with 
Cartan decomposition $\fra g = \fra k \oplus \fra p$. For $1 \leq s \leq 3$, 
let $\Cal O_s \subseteq \fra p$ be the subspace of Jordan rank $s$.
For $\fra g = \fra{sp}(3,\Bobb R)$ and $\fra g = \fra{su}(3,3)$,
the Jordan rank amounts to the ordinary rank of a matrix, where 
$\roman S_{\Bobb C}^2[\Bobb C^3]$ is identified with the symmetric 
$(3\times 3)$-matrices. For $\fra g = \fra{so}^*(12)$, when we identify
$\Lambda^2[\Bobb C^6]$ with the skew-symmetric $(6\times 6)$-matrices,
the Jordan rank amounts to one half the ordinary rank of a matrix.
The resulting decomposition of $\fra p$ is a stratification
whose strata coincide with the $K^{\Bobb C}$-orbits,
and the closures constitute an ascending sequence
$$
\{0\}
\subseteq \overline {\Cal O_1}
\subseteq \overline {\Cal O_2}
\subseteq \overline {\Cal O_3} = \fra p
\tag4.1
$$
of complex affine varieties. Projectivization yields the ascending sequence
$$
Q_1 \subseteq 
Q_2 \subseteq  Q_3 =\bold P(\fra p).
\tag4.2
$$
As far as the complex analytic structures are concerned, this is the 
sequence (1.6). By construction, $Q_2$ is a determinantal cubic, the 
requisite determinant over the octonions being that introduced by 
Freudenthal \cite\freudone, and $Q_1$ is the corresponding Severi variety.
In fact, $Q_1$ is the {\it closed\/} $K^{\Bobb C}$-orbit in $\bold P(\fra p)$, 
and the homogeneous space descriptions (1.1)--(1.4) are immediate. For
$\fra g =\fra e_{7(-25)}$, the cubic $Q_2$ is the (projective) generic norm 
hypersurface, sometimes referred to in the literature as 
{\it Freudenthal\/} cubic; it has been studied  already by E. Cartan, though.
Jacobson has shown that this cubic is rational \cite\jacobstw.
In the next section we recall briefly how the two ascending sequences
(4.1) and (4.2) arise.

\medskip\noindent
{\bf 5. Holomorphic nilpotent orbits and the proof of Theorem 1.5}
\smallskip\noindent
A detailed account of holomorphic nilpotent orbits is given
in our paper \cite\kaehler.
Here we recall only what we need for ease of exposition.

Let $(\fra g, z)$ be a semisimple Lie algebra of hermitian type,
with Cartan decomposition $\fra g = \fra k \oplus \fra p$.
We refer to an  adjoint orbit $\Cal O \subseteq \fra g$
having the property that the projection map from $\fra g$ to $\fra p$,
restricted to $\Cal O$, is a diffeomorphism onto its image,
as a {\it pseudoholomorphic\/} orbit. A pseudoholomorphic orbit $\Cal O$
inherits a complex structure from the complex structure $J_z$ on $\fra p$,
and this complex structure, combined with the Kostant-Kirillov-Souriau 
form on $\Cal O$, viewed as a coadjoint orbit by means of 
(a positive multiple of) the Killing form, turns
$\Cal O$ into a (not necessarily positive) K\"ahler manifold.

We now choose a positive multiple of the Killing form. We say that
a pseudoholomorphic orbit $\Cal O$ is {\it holomorphic\/} provided
the resulting K\"ahler structure on $\Cal O$ is positive. The name 
\lq\lq holomorphic\rq\rq\ is intended to hint at the fact that
the holomorphic discrete series representations of $G$ arise from holomorphic 
quantization on integral {\it semisimple\/} holomorphic orbits but, beware, 
the requisite complex structure (needed for the construction of the 
holomorphic discrete series representation) is not the one arising from 
projection to $\fra p$.

Let $\Cal O$ be a holomorphic {\it nilpotent\/} orbit, and let 
$C^{\infty}(\overline{\Cal O})$ be the algebra of {\it Whitney smooth 
functions\/} on the ordinary topological closure $\overline{\Cal O}$ of 
$\Cal O$ (not the Zariski closure) resulting from the embedding of 
$\overline{\Cal O}$ into $\fra g^*$. The Lie bracket on $\fra g$ passes to 
a Poisson bracket $\{\cdot,\cdot\}$ on $C^{\infty}(\overline{\Cal O})$,
even though $\overline{\Cal O}$ is {\it not\/} a smooth manifold. This 
Poisson bracket turns $\overline{\Cal O}$ into a stratified symplectic space.

A {\it complex analytic stratified K\"ahler structure\/} on a stratified 
space $N$ is a stratified symplectic structure 
$(C^{\infty}(N),\{\cdot,\cdot\})$ ($N$ is not necessarily smooth and
$C^{\infty}(N)$ not necessarily an algebra of ordinary smooth functions)
together with a complex analytic structure which, on each stratum,
\lq\lq combines\rq\rq\  with the symplectic structure on that stratum
to a K\"ahler structure, in the following sense: The stratification 
underlying the stratified symplectic structure is a refinement of the 
complex analytic stratification whence each stratum is a complex manifold; 
each holomorphic function is a smooth function in $C^{\infty}(N,\Bobb C)$;
and on each stratum, the Poisson  structure is symplectic in such a way taht
that the symplectic structure combines with the complex structure to a 
K\"ahler structure. See \cite\kaehler\ for details. The structure may be 
described in terms of a {\it stratified K\"ahler polarization\/} \cite\kaehler\ 
which induces the K\"ahler polarizations on the strata and encapsulates 
the mutual positions of these polarizations on the strata. A complex 
polarization can no longer be thought of as being given by the 
$(0,1)$-vectors of a complex structure, though. When the complex analytic 
structure is normal we simply refer to a {\it normal K\"ahler structure\/}.
We now recall a few facts from  \cite\kaehler.

\noindent
{(5.1)} {\sl Given a holomorphic nilpotent orbit $\Cal O$, the closure
$\overline{\Cal O}$ is a union of finitely many holomorphic nilpotent orbits.
Moreover, the diffeomorphism from $\Cal O$ onto its image in $\fra p$
extends to a homeomorphism from the closure $\overline{\Cal O}$
onto its image in $\fra p$, this homeomorphism turns $\overline{\Cal O}$
into a complex affine variety, and the complex analytic structure, in turn,
combines with the  Poisson structure 
$(C^{\infty}(\overline{\Cal O}),\{\cdot,\cdot\})$
to a normal (complex analytic stratified) K\"ahler structure.\/}
See Theorem 3.2.1 in \cite\kaehler.

Let  $r$ be the real rank of $\fra g$. There are $r+1$ holomorphic nilpotent 
orbits $\Cal O_0,\dots, \Cal O_r$, and these are linearly ordered
in such a way that
$$
\{0\}=\Cal O_0 \subseteq \overline{\Cal O_1} 
\subseteq \ldots \subseteq \overline{\Cal O_r} = \fra p,
\tag5.2
$$
cf. (3.3.10) in \cite\kaehler. Recall that the Cartan decomposition induces the 
decomposition
$\fra g^{\Bobb C} = \fra k^{\Bobb C}\oplus \fra p^+\oplus \fra p^-$
of the complexification $\fra g^{\Bobb C}$ of $\fra g$,
$\fra p^+$ and $\fra p^-$ being the holomorphic and antiholomorphic 
constituents, respectively, of $\fra p^{\Bobb C}$.

\noindent
{(5.3)}
{\sl The projection from $\overline {\Cal O_r}$ to $\fra p$ is a 
homeomorphism onto $\fra p$, and the $G$-orbit stratification of 
$\overline {\Cal O_r}$ passes to the $K^{\Bobb C}$-orbit stratification
of $\fra p\cong \fra p^+$. Thus, for $1 \leq s \leq r$, restricted to 
$\Cal O_s$, this homeomorphism is a $K$-equivariant diffeomorphism
from $\Cal O_s$ onto its image in $\fra p^+$, and this image is a
$K^{\Bobb C}$-orbit in $\fra p^+$.\/} See Theorem 3.3.11 in \cite\kaehler.

\noindent
{\smc Remark 5.4.} The holomorphic nilpotent orbits in 
a simple Lie algebra $\fra g$ of hermitian type are precisely those which 
have the property that the projection to $\fra p\cong\fra p^+$ realizes the
Kostant-Sekiguchi correspondence. See Remark 3.3.13 in \cite\kaehler.  

\demo{Proof of Theorem 1.5}
The ascending sequence (5.3)  of affine complex varieties 
determines the ascending sequence
$$
Q_1 \subseteq Q_2 \subseteq  \dots \subseteq Q_r=\bold P(\fra p)
\tag5.5
$$
of projective varieties where $\bold P(\fra p)$ is the projective space
on $\fra p$ and where each $Q_s$ ($1 \leq s \leq r$) arises from 
$\overline{\Cal O_s}$ by projectivization. The stratified K\"ahler structures 
on the constituents of (5.3) pass to normal K\"ahler structures on the 
constituents of (5.5). See Section 10 in \cite\kaehler\ for details.

Exploiting the sequence (5.5) for each of the four simple rank 3 
hermitian Lie algebras, from  the description of Severi varieties 
in terms of the simple regular rank 3  hermitian Lie algebras
given in Section 4, we obtain a proof of Theorem 1.5 in the 
introduction. In particular, the discussion in Section 8 (Appendix) of 
\cite\atiybern\ reveals that, under the circumstances of Theorem 1.5,
$Q_2$ is indeed the chordal variety of the Severi variety $Q_1$. \qed
\enddemo

In each of these cases, the cubic polynomial which determines the 
hypersurface $Q_2$ is the fundamental relative invariant 
(Bernstein-Sato polynomial) of the corresponding (irreducible regular) 
pre-homogeneous space, cf. Theorem 7.1 in \cite\kaehler\ for the cases 
(1.1)--(1.3) and Theorem 8.4.1 in \cite\kaehler\ for the case (1.4); in the 
latter case, the cubic polynomial is the generic norm or Freudenthal's 
generalized determinant mentioned in Section 4.

\medskip\noindent {\bf 6. GIT-quotients and K\"ahler reduction}
\smallskip\noindent
Given a Hodge manifold $N$, endowed with an appropriate group of symmetries 
and momentum mapping, reduction carries it to a complex analytic stratified 
K\"ahler space $N^{\roman{red}}$ which is as well a projective variety 
\cite\kaehler\ (Section 4) and the question arises whether a {\sl complex 
projective space into which $N^{\roman{red}}$ embeds carries an exotic 
structure which, via the (Kodaira) embedding, restricts to the 
complex analytic stratified K\"ahler structure on\/} $N^{\roman{red}}$.
We now describe a situation where the answer to this question is positive
and, furthermore, we explain how, for the three classical cases,
the circumstances of Theorem 1.5 may be subsumed under this situation.
More details and more general results may be found in \cite\kaehler.

Let $\Bobb K = \Bobb R, \Bobb C, \Bobb H$, and consider the standard 
(right) $\Bobb K$-vector space $\Bobb K^s$ of dimension $s$,
endowed with a (non-degenerate) positive definite hermitian form
$(\cdot,\cdot)$; further, let $V = \Bobb K^n$,  endowed with a skew form 
$\Cal B$ and compatible complex structure $J_V$  such that associating
$\Cal B(u,J_V v)$  to $u,v \in V$ yields a positive definite hermitian form 
on $V$. Furthermore, let $\GHH= \GHH(s)= \roman U(V^s,(\cdot,\cdot))$,
$\fra h = \roman{Lie}(\GHH)$, $\HG = \roman U(V,\Cal B)$,
and $\fra g = \roman{Lie}(\HG)$, and denote the split rank of 
$\HG = \roman U(V,\Cal B)$ by $r$. More explicitly:

\noindent
(6.1) $\Bobb K=\Bobb R,\ V = \Bobb R^n,\
n = 2 \ell,\ 
\GHH= \roman O(s),\ 
\HG= \roman {Sp}(\ell,\Bobb R),\ r = \ell$; 
\newline\noindent
(6.2) $\Bobb K=\Bobb C,\ V = \Bobb C^n,\ 
n = p +q,\ 
\GHH= \roman U(s),\ 
\HG= \roman U(p,q)$, $p \geq q$, $r = q$;
\newline\noindent
(6.3) $\Bobb K=\Bobb H,\ V = \Bobb H^n,\ 
\GHH= \roman U(s,\Bobb H) = \roman {Sp}(s),\ 
\HG= \roman O^*(2 n)$, $r = [\frac n2]$.

Let $W=W(s)=\roman{Hom}_{\Bobb K}(\Bobb K^s,V)$. The 2-forms
$(\cdot,\cdot)$ and $\Cal B$ induce a symplectic structure $\omega_W$ on $W$,
and the groups $\KK$ and $\HG$ act on $W$ in an obvious fashion:
Given $x \in \KK, \ \alpha \in \roman{Hom}_{\Bobb K}(\Bobb K^s,V),\ y \in \HG$,
the action is given by the assignment to $(x,y,\alpha)$ of $y \alpha x^{-1}$.
These actions preserve the symplectic structure $\omega_W$ and are hamiltonian:
Given $\alpha \colon \Bobb K^s \to V$, define
$\alpha^\dagger \colon V \to \Bobb K^s$ by
$$
(\alpha^\dagger \bold u,\bold v) = \Cal B(\bold u,\alpha \bold v),
\ \bold u \in V,\, \bold v \in \Bobb K^s.
$$
Under the identifications of $\fra {\kk}$ and $\fra {\hg}$ with their duals 
by means of the half-trace pairings, the momentum mappings $\mu_{\KK}$ and
$\mu_{\HG}$ for the $\KK$- and $\HG$-actions on $W$ are given by
$$
\align
\mu_{\KK} \colon W &@>>> \fra {\kk},
\quad
\mu_{\KK}(\alpha)= -\alpha^\dagger \alpha \colon \Bobb K^s \to \Bobb K^s,
\\
\mu_{\HG} \colon W &@>>> \fra {\hg},
\quad
\mu_{\HG}(\alpha) =
 \alpha \alpha^\dagger \colon V \to V
\endalign
$$
respectively; these are the  Hilbert maps of invariant theory for the 
$\KK$- and $\HG$-actions. Moreover, the groups $G$ and $H$ constitute 
a reductive dual pair in $\roman{Sp}(W,\omega_W)$, $J_V$ is a member of \
$\fra {\hg}$, and $(\fra g,\frac 12 J_V)$ is a simple (or reductive with 
simple semisimple constituent) Lie algebra of hermitian type. For 
$\phi \colon \Bobb K ^s \to V$, let $J_W (\phi) = J_V \circ \phi$.
This yields a complex structure $J_W$ on $W$ which,
together with the symplectic structure $\omega_W$, turns $W$ into a flat
(not necessarily positive) K\"ahler manifold.
When we wish to emphasize in notation that $W$ is endowed
with this complex structure we write $W_J$.
The $H$-action on $W_J$ preserves the complex structure.
See Section 5 of \cite\kaehler\ for details.

Since the $G$- and $H$-actions centralize each other, the $G$-momentum mapping
$\mu_{\HG}$ induces a map $\overline \mu_{\HG}$ from the ${\KKK}$-reduced space
$W(s)^{\roman{red}} = \mu_{\KKK}^{-1}(0)\big / {\KKK}$ into $\fra {\hg}$.  
We now recall the following facts from Theorem 5.3.3 in \cite\kaehler.

\noindent (6.4)
{\rm (1)} 
{\sl The induced map $\overline \mu_{\HG}$ from the ${\KKK}$-reduced space
$W(s)^{\roman{red}}$ into $\fra {\hg}$  is a proper embedding of
$W(s)^{\roman{red}}$ into $\fra {\hg}$ which, for $s \leq r$,
induces an isomorphism of normal K\"ahler spaces
from $W(s)^{\roman{red}}$ onto the closure $\overline {\Cal O_{s}}$
of the holomorphic nilpotent orbit\/} $\Cal O_{s}$.
\newline\noindent
{\rm (2)} {\sl For $2\leq s \leq r$, the injection of $W(s-1)$ into $W(s)$ 
induces an injection of $W(s-1)^{\roman{red}}$ into $W(s)^{\roman{red}}$
which, under the identifications with the closures of holomorphic
nilpotent orbits, amounts to the inclusion\/}
$\overline {\Cal O_{s-1}} \subseteq \overline {\Cal O_s}$.
\newline\noindent
{\rm (3)} {\sl For $s > r$,
the obvious injection of $W(s-1)$ into $W(s)$ induces an isomorphism
of $W(s-1)^{\roman{red}}$ onto $W(s)^{\roman{red}}$.\/}
\newline\noindent
{\rm (4)} {\sl For $s \geq 1$, under the projection from 
$\fra g =\fra k \oplus \fra p$ to $\fra p$, followed by the identification 
of the latter with $\fra p^+$, the image of the reduced space
$W^{\roman{red}}$ in $\fra g$ is identified with the affine complex 
categorical quotient $W_J\big /\big / {\KKK}^{\Bobb C}$, realized in
the complex vector space $\fra p^+$.\/}

The projective space $\Bobb P[W]$ being endowed with a positive multiple of 
the Fubini-Study metric, the momentum mapping passes to a momentum mapping
$\mu_{\KK} \colon \Bobb P[W] \to \fra {\kk}$. For $s=r$,
K\"ahler reduction yields the projective space $\Bobb P[\fra p]$ on $\fra p$,
endowed with an exotic K\"ahler structure; for $1 \leq s \leq r$, we obtain 
a normal K\"ahler space $Q_s$ together with an ascending sequence
$$
Q_1 \subseteq Q_2 \subseteq \ldots\subseteq  Q_s
$$
of normal K\"ahler spaces which are, in fact, the closures of the strata of 
$Q_s$; (6.4) above implies that this sequence amounts to the sequence (5.5),
truncated at $Q_s$. Moreover, the embeddings of the $Q_s$'s into 
$Q_r=\Bobb P[\fra p]$ are Kodaira embeddings. In view of (6.4) above,
the normal K\"ahler structure coming from the K\"ahler reduction procedure
coincides with the structure coming from projectivization of
the closure of the corresponding holomorphic nilpotent orbit.
See Section 10 of \cite\kaehler\ for details.

For $r=3$, this construction recovers the sequence (1.6). For intelligibility, 
we describe briefly the various constituents where the numbering
(6.5.*) corresponds to the numbering (1.*) for $*=1,2,3$:
\newline\noindent
The sequence $\Bobb P[W(1)] \subseteq \Bobb P[W(2)] \subseteq \Bobb P[W]$ has 
the form:
\newline\noindent
(6.5.1) 
$\Bobb P^2 \Bobb C \subseteq \Bobb P^5 \Bobb C \subseteq \Bobb P^8 \Bobb C$
\newline\noindent
(6.5.2) 
$\Bobb P^5 \Bobb C \subseteq \Bobb P^{11} \Bobb C \subseteq\Bobb P^{17}\Bobb C$
\newline\noindent
(6.5.3) 
$\Bobb P^{11}\Bobb C\subseteq\Bobb P^{23}\Bobb C\subseteq\Bobb P^{35}\Bobb C$
\newline\noindent
The sequence (1.6) has the form:
\newline\noindent
(6.5.1$'$) 
$X =\Bobb P^2\Bobb C \subseteq Q^4 \subseteq \Bobb P^5 \Bobb C 
=\Bobb P^8 \Bobb C \big/\big/ \roman O(3,\Bobb C)$
\newline\noindent
(6.5.2$'$) 
$X=\Bobb P^2\Bobb C\times\Bobb P^2\Bobb C\subseteq Q^7\subseteq\Bobb P^8\Bobb C
=\Bobb P^{17} \Bobb C \big/\big/ \roman {GL}(3,\Bobb C)$
\newline\noindent
(6.5.3$'$) 
$X=\roman G_2(\Bobb C^6) \subseteq Q^{13}\subseteq\Bobb P^{14}\Bobb C
=\Bobb P^{35} \Bobb C \big/\big/ \roman {Sp}(3,\Bobb C)$
\newline\noindent
As GIT-quotients, the Severi varieties $X$ and the cubics $Q^*$ may 
be written in the form:
\newline\noindent
(6.5.1$''$) 
$X=\Bobb P^2 \Bobb C \big/\big/ \roman {O}(1,\Bobb C) \cong \Bobb P^2\Bobb C$,
$Q^4 =\Bobb P^5 \Bobb C \big/\big/ \roman O(2,\Bobb C)$
\newline\noindent
(6.5.2$''$) 
$X=\Bobb P^2\Bobb C\times\Bobb P^2\Bobb C=
\Bobb P^5 \Bobb C \big/\big/ \roman {GL}(1,\Bobb C)$,
$Q^7 =
\Bobb P^{11} \Bobb C \big/\big/ \roman {GL}(2,\Bobb C)$,
\newline\noindent
(6.5.3$''$) 
$X=\roman G_2(\Bobb C^6) = 
\Bobb P^{11} \Bobb C \big/\big/ \roman {Sp}(1,\Bobb C)$,
$Q^{13} =\Bobb P^{23} \Bobb C \big/\big/ \roman {Sp}(2,\Bobb C)$.
\newline\noindent
Here, for $1 \leq s \leq 3$, the actions of the groups
$\roman O(s,\Bobb C)$, $\roman {GL}(s,\Bobb C)$, $\roman {Sp}(s,\Bobb C)$
on the corresponding projective spaces arise from the actions of the groups 
$H$ on $W(s)$ listed in (6.1)--(6.3) above. The sequences (6.5.1)--(6.5.3) 
may be viewed as resolutions of singularities (in the sense of stratified 
K\"ahler spaces) for the corresponding sequences (1.6). The disjoint union 
$$
\Bobb P[W]= H^{\Bobb C}\Bobb P[W(1)] \cup
H^{\Bobb C}(\Bobb P[W(2)] \setminus \Bobb P[W(1)])\cup
(\Bobb P[W] \setminus H^{\Bobb C}\Bobb P[W(2)])
$$
is the  $H^{\Bobb C}$-orbit type decomposition of $\Bobb P[W]$ in each case;
here $H^{\Bobb C}$ denotes the complexification of $H$.

In particular, the cubic $Q_2$
(written as $Q^4$, $Q^7$, $Q^{13}$ according to the case considered
where the superscript indicates the complex dimension)
arises by K\"ahler reduction, applied to the projective space $\Bobb P[W(2)]$.
This is an instance of the situation referred to at the beginning
of this section. Furthermore, the Severi variety $Q_1$
arises by K\"ahler reduction, applied to
the projective space $\Bobb P[W(1)] =\Bobb P[V]$.

\noindent
{\smc Remark 6.6.}
In the special case where 
$(\fra g, \fra h) = (\fra {sp}(\ell,\Bobb R), \fra {so}(s,\Bobb R))$,
the space $W$ may be viewed as the unreduced phase space of $\ell$ 
particles in $\Bobb R^s$, and $\mu_{\GH}$ is the angular momentum mapping. 
Thus the normal K\"ahler space $\overline {\Cal O_s}$ arises as the reduced 
phase space of a system of $\ell$ particles in $\Bobb R^s$ with total 
angular momentum zero. Likewise, the compact normal K\"ahler space $Q_s$ 
which, complex analytically, is a projective variety, arises as the reduced 
phase space of a system of $\ell$ harmonic oscillators in $\Bobb R^s$ with 
total angular momentum zero and constant energy. Here the energy is encoded 
in the stratified symplectic Poisson structure; changing the energy amounts 
to rescaling the Poisson structure. The constituents $Q_s$ 
($1 \leq s \leq \ell$) of the ascending sequence (5.5) have the following 
interpretation: The top stratum $Q_{\ell} \setminus Q_{\ell-1}$ consists of 
configurations in general position, that is, $\ell$ harmonic oscillators in 
$\Bobb R^{\ell}$ with total angular momentum zero
such that the positions and momenta
do not lie in a cotangent bundle $\roman T^* \Bobb R^s$ for some 
$\Bobb R^s \subseteq \Bobb R^{\ell}$ with $s <\ell$. For $1 \leq s<\ell$, 
the stratum $Q_s \setminus Q_{s-1}$ consists of configurations of $\ell$ 
harmonic oscillators in $\Bobb R^{\ell}$ 
with total angular momentum zero such that the 
positions and momenta lie in the cotangent bundle $\roman T^* \Bobb R^s$
for some $\Bobb R^s \subseteq \Bobb R^{\ell}$ but {\it not\/} in a 
cotangent bundle of the kind $\roman T^* \Bobb R^{s-1}$, whatever 
$\Bobb R^{s-1} \subseteq \Bobb R^{\ell}$; here the convention is that 
$Q_0$ is empty. Under these circumstances, the complex analytic structure
does not have an immediate mechanical significance but helps understanding 
the geometry of the reduced phase space. The complex analytic structure is 
important for issues related with quantization, cf. \cite\kaehredu.
Though this is not relevant here, we note that,
for $s>\ell$, the obvious map from $Q_{\ell}$ to $Q_s$ is an isomorphism of 
stratified K\"ahler spaces, and no new geometrical phenomenon occurs.
Since the other dual pairs lie in some $\fra{sp}(n,\Bobb R)$,
in the other cases (corresponding to (1.2) and (1.3)),
the constituents $Q_s$ ($1 \leq s \leq r$) of the ascending sequence
(5.5) admit as well interpretations in terms of suitable constrained 
mechanical systems.

The exceptional Severi variety (1.4) is notably absent. The question 
whether this variety and the corresponding ambient spaces $Q_2$ and 
$Q_3$  arise from K\"ahler reduction in the same way as the varieties 
(1.1)--(1.3) and corresponding ambient spaces is related with that of existence
of a dual pair $(G,H)$ with $G=E_{7(-25)}$ and $H$ compact but apparently 
there is no such dual pair. It would be extremely interesting to develop an 
alternative construction which yields the ascending sequence (1.6) 
for the exceptional Severi variety by K\"ahler reduction
applied to a suitable K\"ahler manifold, perhaps more complicated
than just complex projective space.

\medskip
\centerline{\smc References}
\medskip
\widestnumber\key{999}

\ref \no \atiybern
\by M. Atiyah and J. Berndt
\paper Projective planes, Severi varieties and spheres
\linebreak
\finalinfo{\tt math.DG/0206135}
\endref

\ref \no \farakora
\by J. Faraut and A. Koranyi
\book Analysis On Symmetric Cones
\bookinfo Oxford Mathematical Monographs
 No. 114
\publ Oxford University Press
\publaddr Oxford U.K.
\yr 1994
\endref

\ref \no \freudone
\by H. Freudenthal
\paper Beziehungen der $E_7$ und $E_8$ zur Oktavenebene
\jour Indagationes Math.
\vol 16
\yr 1954
\pages 218--230
\endref

\ref \no \hartsboo
\by  R. Hartshorne
\book Algebraic Geometry
\bookinfo Graduate texts in Mathematics
 No. 52
\publ Springer
\publaddr Berlin-G\"ottingen-Heidelberg
\yr 1977
\endref

\ref \no \howeone
\by R. Howe
\paper Remarks on classical invariant theory
\jour  Trans. Amer. Math. Soc.
\vol 313
\yr 1989
\pages  539--570
\endref

\ref \no \poiscoho
\by J. Huebschmann
\paper Poisson cohomology and quantization
\jour 
J. f\"ur die reine und angewandte Mathematik
\vol  408 
\yr 1990
\pages 57--113
\endref
\ref \no  \souriau
\by J. Huebschmann
\paper On the quantization of Poisson algebras
\paperinfo Symplectic Geometry and Mathematical Physics,
Actes du colloque en l'honneur de Jean-Marie Souriau,
P. Donato, C. Duval, J. Elhadad, G.M. Tuynman, eds.;
Progress in Mathematics, Vol. 99
\publ Birkh\"auser Verlag
\publaddr Boston $\cdot$ Basel $\cdot$ Berlin
\yr 1991
\pages 204--233
\endref

\ref \no  \oberwork
\by J. Huebschmann
\paper Singularities and Poisson geometry of certain representation spaces
\paperinfo in: Quantization of Singular Symplectic Quotients,
N. P. Landsman, M. Pflaum, M. Schlichenmaier, eds.,
Workshop, Oberwolfach,
August 1999,
Progress in Mathematics, Vol. 198
\publ Birkh\"auser Verlag
\publaddr Boston $\cdot$ Basel $\cdot$ Berlin
\yr 2001
\pages 119--135
\finalinfo{\tt math.DG/0012184}
\endref

\ref \no \kaehler
\by J. Huebschmann
\paper K\"ahler spaces, nilpotent orbits, and singular reduction
\jour Memoirs AMS (to appear)
\finalinfo {\tt math.DG/0104213}
\endref

\ref \no \kaehredu
\by J. Huebschmann
\paper K\"ahler quantization and reduction
\paperinfo {\tt math.SG/0207166}
\endref

\ref \no \lradq
\by J. Huebschmann
\paper Lie-Rinehart algebras, descent, and quantization
\jour Fields Institute Communications (to appear)
\finalinfo {\tt math.SG/0303016}
\endref

\ref \no \jacobstw
\by N. Jacobson
\paper Some projective varieties defined by Jordan algebras
\jour J. of Algebra
\vol 97
\yr 1985
\pages 556--598
\endref

\ref \no \kostasix
\by B. Kostant
\paper The principal three-dimensional subgroup
and the Betti numbers of a complex simple Lie group
\jour Amer. J. of Math.
\vol 81
\yr 1959
\pages 973--1032
\endref

\ref \no \lanmaone
\by J. M. Landsberg and L. Manivel
\paper The projective geometry of Freudenthal's magic square
\jour J. of Algebra
\vol 239
\yr 2001
\pages 477--512
\endref

\ref \no \mccrione
\by  K. McCrimmon
\paper Jordan algebras and their applications
\jour Bull. Amer. Math. Soc.
\vol 84
\yr 1978
\pages 612--627
\endref

\ref \no \satakboo
\by I. Satake
\book Algebraic structures of symmetric domains
\bookinfo Publications of the Math. Soc. of Japan, vol. 14
\publ Princeton University Press
\publaddr Princeton, NJ
\yr 1980
\endref

\ref \no \sekiguch
\by J. Sekiguchi
\paper Remarks on real nilpotent orbits of a symmetric pair
\jour J. Math. Soc. Japan
\vol 39
\yr 1987
\pages 127--138
\endref

\ref \no \severone
\by F. Severi
\paper Intorni ai punti doppi impropi di una superficie generale
dello spazio a quattro dimenzioni, e a'suoi punti tripli apparenti
\jour Rendiconti del Circolo Matematico di Palermo
\vol 15
\yr 1901
\pages 33--51
\endref

\ref \no \zakone
\by F. L. Zak
\paper Severi varieties
\jour Math. USSR, Sb.
\vol 54
\yr 1986
\pages 113--127
\endref

\enddocument